\documentclass[12pt]{article}          
\usepackage{amssymb}

\usepackage{amssymb, amscd}
\usepackage{latexsym}
\usepackage{amsmath}
\usepackage{amssymb}
\usepackage{ifthen}
\usepackage{amsthm}
\newtheorem{dfn}{Definition}[section]

\newtheorem{rem}[dfn]{Remark}
\newtheorem{rems}[dfn]{Remarks}

\newtheorem{lem}[dfn]{Lemma}
\newtheorem{sublem}[dfn]{Sublemma}
\newtheorem{claim}[dfn]{Claim}
\newtheorem{prop}[dfn]{Proposition}
\newtheorem{cor}[dfn]{Corollary}

\newtheorem{ques}[dfn]{Question}

\def\proof{\par\medskip\noindent{\it Proof: }}

\def\A{{\cal A}}

\def\H{{\Bbb H}}
\def\Z{{\Bbb Z}}

\def\E{{\cal E}}

\def\T{{{\cal T}}}
\def\U{{{\cal U}}}
\def\P{{{\cal P}}}

\def\al{\alpha}
\def\be{\beta}
\def\ga{\gamma}
\def\g{\gamma}

\def\G{\Gamma}

\def\De{\Delta}
\def\Si{\Sigma}
\def\si{\sigma}

\def\D{{\cal D}}

\def\La{\Lambda}

\def\si{\sigma}
\def\ve{\varepsilon}
\def\sie{\Si_{\ve}}
\def\gsi{G_{\si}}
\def\r{\rho}

\def\D{\partial}

\def\hook{\hookrightarrow}

\def\iso{{\rm Isom}({\H^n})}
\def\is3{{\rm Isom}({\H^3})}
\def\hn{{\H^n}}
\def\h3{{\H^3}}
\def\h3{{\H^3}}
\def\oh3{\overline{\H^3}}

\def\ag{A_g}

\def\Mth{M_{\mu{\rm thin}}}
\def\Mtk{M_{\mu{\rm thick}}}
\def\mtk{M_{{\mu\over 2} thick}}
\def\Vo{\rm Vol}
\def\ti{\widetilde}
\def\g+{ g^{+}}
\def\pim{\pi_1(M)}

\def\pims{\pi_1(M)'}
\def\pimtk{\pi_1\Mtk}

\def\qed{\ {\it QED.}}
\def\hook{\hookrightarrow}
\def\bx{$\hfill\blacksquare$}
\begin{document}

\title{On the Complexity and Volume of Hyperbolic 3-Manifolds.}
\author{Thomas Delzant and Leonid Potyagailo}
\date{}
\maketitle \noindent

\bigskip

\begin{abstract}

\noindent We compare the volume of a  hyperbolic 3-manifold $M$
of finite volume    and the  complexity  of its fundamental group. 
\end{abstract}
\bigskip
\bigskip
\section{Introduction.}

\noindent {\bf Complexity  of 3-manifolds and groups.} One of the
most striking corollaries of the recent solution of the
geometrization conjecture for 3-manifolds is the fact that every
aspherical 3-manifold is uniquely determined by its fundamental
group.  It seems to be natural to think that a
topological/geometrical description of a 3-manifold $M$ produces
the simplest way to describe its fundamental group $\pi_1(M);$ on
the other hand, the simplest way to define the group $\pi_1(M)$
gives rise to the most efficient way to describe $M.$ More
precisely, we want to  compare the complexity of 3-manifolds and
their fundamental groups.

The study of the
 complexity   of 3-manifolds goes back
  to the classical work
 of   H. Kneser \cite{K}. Recall that the Kneser complexity  invariant
 $k(M)$ is defined to be the minimal number   of simplices of a
 triangulation of the manifold $M$. The main result of Kneser is that
  this complexity serves as a bound of the number of embedded incompressible
  2-spheres in $M$, and bounds the numbers of factors in a decomposition of $M$
  as a connected sum. A version of
  this complexity
 was used by W. Haken to prove  the existence of hierarchies
for a large class of compact 3-manifolds (called since then Haken
manifolds). Another measure of the complexity
  $c(M)$ for the 3-manifold
 $M$ is  due to
 S.~Matveev. It is the minimal number of vertices of a special spine of
 $M$ [Ma]. It is shown that in many important cases (e.g. if $M$ is
  a non-compact hyperbolic 3-manifold of finite volume) one has $k(M)=c(M)$ \cite{Ma}.

 The rank (minimal number of generators) is also a measure of complexity of a
  finitely generated group. According   to the classical theorem of I.~Grushko
   \cite{Gr},
the rank of a free product of groups is the sum of their ranks.
This immediately implies that every finitely generated group is a
free product of finitely many freely indecomposible factors, which
is an algebraic analogue of Kneser theorem.

  For a finitely presented group $G$ a measure of complexity of $G$
   was defined in \cite{De}. Here is its definition :

\begin{dfn}
\label{relinv}

Let $G$ be a finitely presented group. We say that $T(G) \leq t$ if there exists a
simply-connected 2-dimensional complex $P$ such that $G$ acts
freely  and simplicially on $P$ and the
 the number of 2-faces of the quotient  $\Pi=P/G$
  is less than $t.$

  \medskip
If the group $G$ is defined by a presentation $<a_1,...a_r ; R_1,...R_n>$ the sum $\Sigma (\vert R_i\vert -2)$ serves as a natural bound for $T(G)$.

\end{dfn}

\medskip

 \noindent  Note that an inequality between Kneser complexity and this invariant
is obvious.
 Indeed, by contracting a maximal subtree of the  2-dimensional skeleton
 of a triangulation of $M$ one obtains a triangular presentation of the group $\pi_1(M).$
Since every 3-simplex has four 2-faces it follows

$$T(\pi_1(M)) \leq 4 k(M) .$$

In order to compare the complexity of a manifold and that of its
fundamental group, it is enough to find a function $\theta$ such
that $\theta(\pi_1(M))\leq T(\pi_1(M)).$ Note that the existence
of such a function follows from  G. Perelman's solution of the
geometrization conjecture [Pe 1-3]. Indeed  there could exist at
most finitely many different 3-manifolds having the fundamental
groups isomorphic to the same group $G$ (for  irreducible
3-manifolds with boundary this was shown much earlier in
\cite{Swa}). The question which still remains open is to describe
the asymptotic behavior of the function $\theta$.

  Note that for certain
lens spaces the following inequality is proven in \cite{PP}:

$$c(L_{n,1}) \leq \ln n\approx {\rm const}\cdot T(\Z/n\Z).$$

 However, the above problem remains widely open
for irreducible 3-manifolds with infinite fundamental group. If
$M$ is a compact hyperbolic 3-manifold, D.~Cooper showed \cite{C}:

$${\rm Vol} M \leq \pi\cdot T(\pi_1(M))\hfill\eqno(C).$$

\noindent where ${\rm Vol} M$ is the hyperbolic volume of $M.$
Note that  the converse inequality in dimension 3 is not true:
there exists infinite sequences of different hyperbolic
3-manifolds $M_n$ obtained by Dehn filling on a fixed finite
volume hyperbolic manifold $M$ with cusps such that ${\rm Vol} M_n
< {\rm Vol} M$ \cite{Th}. The ranks of the groups $\pi_1(M_n)$ are
all bounded by ${\rm rank} (\pim)$ and since $\pi_1(M_n)$ are not
isomorphic, we must have $T(\pi_1(M_n))\to\infty.$ So the
invariant $T(\pi_1(M))$ is not comparable with  the volume of
hyperbolic 3-manifolds. This difficulty can be overcome using the
following relative version of the invariant $T$ introduced in
\cite{De}:

\begin{dfn}
\label{relinv}  Let $G$ be a finitely presented group, and $\cal
E$ be a family of subgroups. We say that $T(G, {\cal E}) \leq  t$
if there exists a simply-connected 2-dimensional complex $P$ such
that $G$ acts
 simplicially on $P$, the number of 2-faces of the quotient (an  orbihedron)
  $\Pi=P/G$
  is less than $t$, and the stabilizers of vertices of $P$ are elements of
   $\cal E$.

\end{dfn}
  \medskip

The main goal of the present paper is to obtain uniform constants
comparing the volume   of a  hyperbolic 3-manifold  $M$ of finite
volume and the relative invariant $T(\pi_1 (M), E)$ where $E$ is
the family of its elementary subgroups.

To finish our historical discussion  let us point out that the
relative invariant $T(G, E)$ allows one to prove  the
accessibility of a finitely presented group $G$ without 2-torsion
over elementary subgroups  \cite{DePo1}.
 Using these methods it was shown recently that for
hyperbolic groups without 2-torsion any canonical hierarchy over
finite subgroups and one-ended subgroups is finite \cite{Va}. The
relative invariant $T$ and the hierarchical accessibility was used
in \cite{DePo2} to give  a criterion of the co-Hopf property  for
geometrically finite discrete subgroups of $\iso$.

\bigskip

\noindent {\bf Main Results.} Let $M$ be a hyperbolic 3-manifold
of finite volume. We consider the  family $E_{\mu}$ of all
  elementary subgroups of $\pim$ having translation length less
  than the Margulis constant $\mu=\mu(3)$. The family $E_\mu$
   includes all parabolic
  subgroups of $G$ as well as cyclic loxodromic ones representing
  geodesics in $M$ of length less than $\mu$ (see also
  the next Section).

     The  first result of the
  paper is the following:

\bigskip

\noindent {\bf Theorem A.} {\it There exists a constant $C$
  such that  for every hyperbolic 3-manifold $M$ of finite volume
  the following inequality holds:}
 $$C^{-1} T(G, E_{\mu})  \leq  {\rm Vol} (M) \leq  C T(G, E_{\mu})\hfill \eqno(*)$$

\bx

\medskip

\noindent The following are  corollaries  of  Theorem A.

\medskip

\begin{cor}
\label{cov}

Suppose $\displaystyle M_n\stackrel {f_n}{\longrightarrow} M$ is a
sequence of finite coverings over a finite volume 3-manifold $M$
such that ${\rm deg} f_n\to +\infty.$ Then $T(\pi_1(M_n), E_n)\to
+\infty$, where $E_n$ is the above system of  elementary subgroups
of $\pi_1(M_n)$ whose translation length is less than $\mu.$\bx
\end{cor}

\proof The statement follows immediately from the right-hand side
of (*)
 since ${\rm Vol}
(M_n)\to\infty.$\qed

\medskip

\begin{cor}
\label{surg} Let $M_n$ be a sequence of different hyperbolic
3-manifolds obtained by Dehn surgery on a cusped hyperbolic
3-manifold of finite volume $M$. Then

$$T(\pi_1(M_n), E_n) \leq  C\cdot {\rm Vol} (M)\ < +\infty.$$
\end{cor}

\bx

\proof The left-hand side of (*) gives

$$T(\pi_1(M_n), E_n) \leq   C \cdot {\rm Vol} (M_n),$$

\noindent and by \cite{Th} one has ${\rm Vol}(M_n) < {\rm Vol}
(M).$ \qed

\medskip

As it is  pointed out   in   Corollary \ref{cov} above we must
have $T(\pi_1(M_n))\to +\infty$ for the absolute invariant. Our
next result is the following :

\bigskip

\noindent {\bf Theorem B.} {\it (Generalized Cooper inequality)
Let $E$ be  the   family of elementary subgroups of $G$, then one
has

$${\rm Vol} (M) \leq \pi\cdot  T(\pim, E)\hfill \eqno(**)$$

} \bx

\bigskip

 \noindent Note that Theorem B gives a generalization of the
Cooper inequality (C) for   the relative invariant $T(G, E)$.
Furthermore, if one puts   $E=E_\mu$, then Theorem B implies the
right-hand side of (*) in Theorem A. Theorems A and B together
have several immediate consequences:

\begin{cor}
\label{compar}

For the constant $C$ from Theorem A the following statements hold:
\begin{itemize}

  \item[\sf i)] Let $M$ be a finite volume hyperbolic 3-manifold
  and $E_{\mu}$ and $E$ be the above families of elementary subgroups
of $\pi_1(M).$ Then

$$T(\pi_1(M), E_{\mu}) \leq C\cdot\pi\cdot T(\pi_1(M), E).$$

\medskip

 \item[\sf ii)] Let $M$ be a hyperbolic 3-manifold such that $M=\Mtk$, i.e.
  every loop  in $M$ of   length less than $\mu$ is homotopically trivial. Then

  $$T(\pi_1(M)) \leq C\cdot \pi\cdot T(\pi_1(M), E).$$

\end{itemize}
\end{cor}

\proof i)   By    Theorems  A and B we have

$$T(\pi_1(M), E_{\mu}) \leq C {\rm Vol} (M) \leq C\cdot\pi\cdot T(\pi_1(M), E).\ \
\qed$$

ii) Since $E_{\mu}=\emptyset$ the result follows from i). \qed

\bigskip

Let us now briefly describe the content of the paper. In Section 2
  we provide some preliminary results needed in the
future.  The proof of Theorem B is  given in Section 3, it
provides  a "simplicial blow-up'' procedure for an orbihedron. In
Section 4 we prove the left-hand side of the inequality (*) using
some standard techniques and the results of Section 2. In the last
Section 5 we discuss some open questions related to the present
paper.

\bigskip

 \noindent {\bf Acknowledgements.} During this work both
authors were partially supported by the ANR grant ${\rm
BLAN}~07-2183619.$ The second author is  grateful to Daryl Cooper
for helpful remarks and to the Max-Planck Institute f\"ur
Mathematik in Bonn, where a part of the work was done. Both
authors are thankful to Anna Lenzhen for  corrections  improving
the paper.

\section{Preliminary results.}

\noindent
 \medskip

 \noindent Let us recall few standard  definitions which we will use in the future.
We say that $G$
 {\it splits } as a graph of groups $X_*=(X, (C_e)_{e\in X^1},
(G_v)_{v\in X^0})$ (where $C_e$ and $G_v$ denote respectively edge
and vertex groups of the graph $X$)    if $G$ is isomorphic to the
fundamental group $\pi_1(X_*)$ in the sense of Serre \cite{Se}.
The Bass-Serre tree $T$ is  the universal cover of the graph
$X=T/G$. When $X$ has only one edge, we will say that $G$ splits
as an amalgamated free product (resp. an HNN-extension) if $X$ has
two vertices (resp. one vertex).

  \begin{dfn}
 Let $G$ be a group  acting on a tree $T$. A subset $H$ of $G$ is
 {\it  elliptic} (resp. {\it
hyperbolic}) in $T$ (and in the graph $T/G$)  if $H$ fixes   a
point in $T$ (resp. does not fix a point in $T$). If $T$ is the
Bass-Serre tree of a splitting of $G$ as a graph of groups, $H$ is
elliptic if and only if it is conjugate into a vertex group of
this graph.

\noindent We say that $G$ splits relatively to a family of
subgroups $ (E_1,...E_n)$, or that the pair $(G,(E_i)_{1\leq i\leq
n})$ splits as a graph of groups, if $G$ splits as a graph of
groups such that all the groups $E_i$  are elliptic in this
splitting. A $(G,(E_i)_{1\leq i\leq n})$-tree is a $G-$tree in
which $E_i$ are   elliptic for all i.\bx

\end{dfn}

  \begin{dfn}

    Suppose $G$ splits as a graph of groups

    $$G=\pi_1(X, C_e, G_v)\hfill\eqno(1)$$
  \noindent relatively to a family of subgroups $E_i\ \{i=1,...,n\}$.

  The decomposition (1) such that all edge groups are non-trivial is called
   {\bf reduced} if every
  vertex group $G_v$ cannot be decomposed relatively to the subgroups
  $E_i\in G_v$ as a graph of groups having one of the subgroups  $C_e$ as
   a vertex group.

 The decomposition (1) is called {\bf rigid} if whenever one has
  a
 $(G, (E_i)_{i\in\{1,...,n\}})$-tree $T^*$
       such that the subgroup $C_e$  contains
  a non-trivial edge stabilizer then $C_e$ acts elliptically on
  $T^*.$\bx

\end{dfn}

\noindent It was shown in \cite{De} that  the sum of relative
$T$-invariants of the vertex groups of a reduced splitting is
 less than or equal to the absolute invariant of $G.$

Recall that  the  Margulis constant $\mu=\mu(n)$ is a number for
which any n-dimensional hyperbolic manifold $M$ can be decomposed
into thick and thin parts : $M=\Mtk\bigsqcup \Mth$ such that the
injectiviry radius at each point of $\Mth$ is less than $\mu/2,$
and
  $\Mtk=M\setminus \Mth.$ By the Margulis Lemma   the
  components of $\Mth$ are either parabolic cusps  or regular
  neighborhoods (tubes) of closed geodesics of $M$ of length less than $\mu.$
  We will denote by  $E=E(\pi_1(M))$ (respectively $E_\mu=E_\mu(\pi_1(M))$) the system of elementary
subgroups   of $\pi_1(M)$ (respectively the systems of subgroups
of $\pimtk$). We will need the following:

\begin{lem}
\label{red}

Let  $H$ be a group admitting the following splitting as a graph
of groups:

$$H=\pi_1(X, C_e, G_v),\hfill\eqno(2)$$

\noindent  where each vertex group $G_v$ is   a lattice in $\iso\
(n>2)$ and $C_e\in E(G_v)\ (n>2)$.

Then  (2) is a reduced and rigid splitting of the couple $(H, \E)$
where $\displaystyle \E=\bigcup_v E(G_v)$.

\end{lem}

\begin{rem} The above Lemma will be further used in a very particular
geometric situation when the group $H$ is the fundamental group of
the double of the thick part $\Mtk$ of $M$ along its boundary.\bx
\end{rem}

 \proof We first claim that it is
enough to prove that every vertex group $G_v$ of the graph $X$
cannot split non-trivially over an elementary subgroup. Indeed, if
it is the case then obviously (2) is reduced. If it is not rigid,
then the couple $(H, \E)$ acts on a simplicial tree $T^*$ such
that one of the groups $C_e$ contains an edge stabiliser $C^*_e$
of $T^*$ and therefore acts hyperbolically on $T^*.$ It follows
that the vertex group $G_v$ containing $C_e$ also acts
hyperbolically on $T^*$ and so is decomposable over elementary
subgroups.

Let us now fix a vertex $v$ and set  $G=G_v$. The Lemma now
follows from the following statement:

\begin{sublem} \cite{Be}
\label{notsplit} Let $G$ be the fundamental group of a Riemannian
manifold $M$ of finite volume of dimension $n > 2$ with pinched
sectional curvature within $[a, b]$ for $a\leq b < 0.$ Then $G$
does not split over a virtually nilpotent group.
\end{sublem}

\proof We provide below a direct proof of this Sublemma in the
case of the constant curvature. Suppose, on the contrary, that

$$G=A*_CB\ \ {\rm or}\ \  G=A*_C, \eqno (3)$$

\noindent where $C$ is an elementary subgroup. Let $\tilde C$ be
the  maximal elementary subgroup containing $C$. The group $\tilde
C$ is virtually abelian and contains a maximal abelian subgroup
$\tilde C_0$ of finite index. We have the following

\begin{claim}

\label{sep}

The group $\tilde C_0$ is separable in $G.$
\end{claim}

\proof \footnote{The argument is due to M.~Kapovich and one of the
authors is thankful for sharing it with him (about 20 years ago).}
Recall that the subgroup $\tilde C_0$ is said {\it separable} if
$\forall g\in G\setminus \tilde C_0$ there exists a subgroup of
finite index $G_0 < G$ such that $\tilde C_0 < G_0$ and $g\not \in
G_0.$ Since $\tilde C_0$ is a maximal abelian subgroup of $G$, and
$g\not\in\tilde C_0,$ it follows that there exists $h\in \tilde
C_0$ such that $\gamma=gh_0g^{-1}h_0^{-1}\not=1.$ The group $G$ is
residually finite, so there exists an epimorphism $\tau : G \to K$
to a finite group $K$ such that $\tau(\gamma)\not=1.$ Since
$\tau(\tilde C_0)$ is abelian, $\tau(\gamma)\not\in \tau(\tilde
C_0)$ and the subgroup $G_0=\tau^{-1}(\tau(\tilde C_0))$ satisfies
our Claim. \qed

Denote $C_0=C\cap \tilde C_0$ (the maximal abelian subgroup of
$C$). We have $\displaystyle \tilde C =\bigcup_{i=1}^m c_i C_0\cup
C_0.$ So by the Claim we can find a subgroup of finite index $G_0$
of $G$ containing $C_0$ such that $c_i\not\in G_0\ (i=1,...,m).$
Then $G_0\cap\tilde C= C_0$ is abelian group and by the Subgroup
Theorem \cite{SW} we have that $G_0$ splits as :

$$ G_0=A_0*_{C'_0}B_0\ \ {\rm or}\ \ G_0=A_0*_{C'_0},\eqno (3')$$

\noindent where $C'_0 < C_0$ is also abelian.  Suppose first that
$G_0=A_0*_{C'_0}B_0$, since $G_0$ is not elementary group, one of
the vertex subgroups of this splitting, say $A_0$ is not
elementary too. Then the map $\varphi : G_0 \to
(cA_0c^{-1})*_{C'_0}B_0,\ \ c\in C'_0,$ such that
$\varphi\vert_{A_0}= cA_0c^{-1}$ and $\varphi\vert_{B_0}={\rm id}$
is an exterior automorphism (as $c$ commutes with every element of
$C_0'$) of infinite order. So the group of the exterior
automorphisms ${\rm Out}(G_0)$ is infinite. This contradicts to
the Mostow rigidity as $G_0$ is still a lattice. In the case of
HNN-extension $G_0=A_0*_{C'_0}= < A_0, t\ \vert\
tC'_0t^{-1}=\psi(C'_0)>$ suppose first that  $t$ does not belong
to the centralizer $Z(C'_0)$ of $C'_0$ in $G_0$. Then we put
$\varphi\vert_{A_0}= cA_0c^{-1}$ for some $c\in C'_0$ such that
$[c, t]\not=1$ and $\varphi(t)=t$. Since $t\not\in Z(C'_0)$ we
obtain again that $\varphi $ is an infinite order exterior
automorphism which is impossible. If, finally, $t\in Z(C'_0)$ then
put $\varphi\vert_{A_0}= id$ and $\varphi(t)=t^2$ and it is easy
to see that $G'_0=\varphi(G_0)$ is a subgroup of index 2 of $G_0$
isomorphic to $G_0$. Then ${\rm Vol}(\H^n/\varphi(G_0)) < +\infty$
and again by Mostow rigidity we must have ${\rm Vol}(\hn/G_0)={\rm
Vol} (\hn/\varphi(G_0)),$ and so $\varphi :G_0\to G_0$ should be
surjective. A contradiction. The Sublemma \ref{notsplit} and Lemma
\ref{red} follow. \qed

\section {Proof of the generalized Cooper inequality.}

 The aim of
this Section is to prove   Theorem B stated in  the Introduction:

\medskip
\noindent {\bf Theorem B.} {\it Let $E$ be an arbitrary family of
  elementary subgroups of $G$, then

$${\rm Vol} (M) \leq \pi\cdot  T(\pim, E)\hfill \eqno(1)$$

} \bx

\proof If  $E=\emptyset,$ then  ${\rm Vol} (M) < \pi\cdot (L-2n),$
where $L$ is the sum of the word-lengths of the relations of
$\pim$ and $n$ is the number of relations  \cite{C}. Let $D$ be a
disk representing a relation  in the presentation complex $R$ of
$\pim$. Then, triangulating $D$ by triangles having vertices on
$\D D,$ we obtain $\vert D\vert-2$ triangles.  So $L-2n$
represents the total number of triangles in $R.$ Thus Cooper's
result implies ${\rm Vol}(M) \leq \pi\cdot T(\pim)$.

Suppose now that $M=\H^3/G$ where $G<\is3$ is a lattice (uniform
or not) and let $E$ be  a   family  of
  elementary subgroups of $G$. Let $P$ be a simply-connected 2-dimensional polyhedron
admitting a simplicial action of $G$ such that the vertex
stabilizers are elements  of    the system $E.$ Let us also assume
that the quotient $\Pi=P/G$ is a finite orbihedron. We will need
the following:

\begin{lem}
\label{emb} There exists a $G$-equivariant simplicial continuous
map $f:P\to \h3\cup
\partial \h3$  such that the images of the
2-simplices of $P$ are   geodesic triangles or ideal triangles of
$\h3$.
\end{lem}

 \proof Let us first construct a $G$-equivariant
continuous  map $f:P\to \overline{\h3}=\h3\cup \partial\h3$ such
that the image of the fixed points for the action $G$ on $P$
belong to $\D\h3.$ To do it we apply the construction from [DePo,
Lemma 1.6] where instead of a tree as the goal space we will use
the hyperbolic space $\h3.$ Let us first construct  a map $\rho
:E\to \h3$ as follows. Since the group $G$ is torsion-free we can
assume that all non-trivial groups in $E$ are infinite. Then for
every elementary group $E_0\in E$ we put $\rho(E_0)=x\in\D\h3$ to
be one of the  fixed points for the action of $E_0$ on $\D\h3$ (by
fixing  a point $O\in\D\h3$ for the image of the trivial group
$\rho(id)).$ The map $\rho$ has the following obvious properties :
\bigskip

\begin{itemize}

  \item[\sf a)] $\forall E_1, E_2\in E\ {\rm if}\ E_1\cap E_2\not=\emptyset\
{\rm then}\ \r(E_1)=\r(E_2)$;

\medskip

  \item[\sf b)] if $\tilde E_0$ is a maximal elementary subgroup then
$\rho(E_0)=\r(\tilde E_0)$ and $\r(g\tilde E_0g^{-1})=g\r(\tilde
E_0)\\ (g\in G).$

\end{itemize}
\bigskip

\noindent We now choose the set of $G$-non-equivalent vertices
$\{p_1,...,p_l\}\subset P$ representing all vertices of $\Pi=P/G$.
 We first construct a map $f$ on zero-skeleton $P^{(0)}$ of the
 complex $P$ by putting $f(p_i)=\rho(E_i)$ and then extend it
 equivariantly $f(gp_i)=gf(p_i) (g\in G).$

 Suppose now $y=(q_1, q_2)\ (q_1,\ q_2\in P^{(0)})$ is an edge of
$P$. To define $f$ on $y$ we distinguish two cases: 1) $H={\rm
Stab} (y)\not=1$ and 2) $H=1.$

\medskip

In the first case  we have necessarily that $E_{g_1}\cap
E_{q_2}=H_0$ is an infinite elementary group where $E_{q_i}$ is
the stabilizer of $q_i.$ Then there exist $g_i\in G$ such that
$q_i=g_i(p_{k_i})\ (i=1,2).$ So $E_{q_i}=g_iE_{p_{k_i}}g_i^{-1}$
and $g_1E_{p_1}g_1^{-1}\cap g_2E_{p_2}g_2^{-1}=H_0.$ It follows
that $E_{p_1}\cap g_1^{-1}g_2E_{p_1}g_2^{-1}g_1$ is an infinite
group and, therefore $f(p_1)=g_1^{-1}g_2(f(p_2))$ implying that
$$f(q_1)=f(g_1p_1)=f(g_2p_2)=f(q_2).$$

\noindent  In the case 2) the stabilizer of the infinite geodesic
 $l=]f(q_1), f(q_2)[\subset \P$ is trivial so we extend $f:y\to l$
 by a piecewise-linear homeomorphism. Having defined  the map $f$ as above on the maximal set of
 non-equivalent edges of $P^{(1)}$ under $G$,  we  extend it
 equivariantly to the 1-skeleton $P^{(1)}$ by putting $f(gy)=gf(y)\ (g\in G).$
Finally we extend $f$
 piecewise linearly  to the 2-skeleton $P^{(2)}.$

We  obtain  a $G$-equivariant continuous  map $f:P\to\oh3$ such
that the all 2-faces of the simplicial complex $f(P)\cap\h3$ are
ideal geodesic triangles.
 The Lemma is proved.  \qed

\begin{rems}
\label{prop}

1. Note that the above Lemma is true in any dimension. We
restricted our consideration to dimension 3 since the further
argument will only concern this case.

\medskip

2. If   the system $E$ contains only parabolic subgroups one can
claim that the action of $G$ on $f(P)\cap\h3$ is in addition
proper. Indeed,  using the convex hull $\P\subset \h3$ of the
maximal family of non-equivalent parabolic points constructed in
\cite{EP} the above argument gives the map
$f:P\to\overline\P\subset \oh3$. By [EP, Proposition 3.5] the set
of faces of $\P$ is locally finite in $\h3.$ Since the boundary of
each face of the 2-orbihedron $f(P)$ constructed above belongs to
$\D\P,$ we obtain that the set of 2-faces of $f(P)\subset\h3$ is
locally finite  in this case.

 \bx

\end{rems}

\noindent If now $W$ is the set of the fixed points for the action
of $G$ on   $P$, we put $P'=P\setminus W$ and
$Q'=f(P')=f(P)\cap\h3.$ Let also  $\nu :P\to \Pi$ and $\pi :\h3\to
M=\h3/G$ denote the natural projections. Then by Lemma \ref{emb}
the map $f$ projects to a simplicial map $F:(\Pi'=P'/G)\to
Q'/G\subset M$ such that the following diagram is commutative:

$$ \begin{CD}\displaystyle   P' @>f\vert_{P'}>>   Q'\subset\h3\\ @V{\nu}VV
@V{\pi}VV
\\  \Pi' @>F>>  Q'/G\subset M \\
\end{CD}
$$

\bigskip

\noindent  Note that, if   $\Pi$ is a simplicial polyhedron, it is
proved in \cite{C} that the hyperbolic area of  $F(\Pi)$   bounds
the volume of the manifold $M.$ This argument does not work if
$\Pi$ is an orbihedron but not a polyhedron. Indeed the complex
$Q'$ above is not necessarily simply connected. So the group $G$
is not isomorphic to $\pi_1(Q'/G)$ but is a non-trivial quotient
of it. Our goal now is to construct a new simplicial polyhedron
$\Sigma$ with the fundamental group $G$ whose image into $M$ has
 area arbitrarily close to that of $F(\Pi')$.  So the main step
in the proof of Theorem B is the following :

\begin{prop}
\label{polappox} (simplicial blow-up procedure). For every $\ve
>0$ there exists a $2$-dimensional complex $\Si_{\ve}$ and a
simplicial map $\varphi_{\ve}:\sie\to M$ such that

\begin{itemize}

\item[1)] The induced map $\varphi_{\ve}:\pi_1\sie\to M$ is  an
isomorphism.

\medskip

and

\medskip

\item[2)] For the hyperbolic area one has:

$$\vert {\rm Area} (\varphi_{\ve}(\sie)) - {\rm Area}(F(\Pi'))\vert
< \ve.$$
\end{itemize}
\end{prop}

\noindent {\it Proof of the Proposition:} Let $\Pi$ be a finite
orbihedron with elementary vertex groups and such that $\pi_1^{\rm
orb}(\Pi)\cong G.$   Let us fix a vertex $\si$ of $\Pi$
 and let $\ti\si\in\nu^{-1}(\si)$ be its  lift in $P.$ We
denote by $\gsi$ the group of the   vertex $\si$ in $G$. By Lemma
\ref{emb} the point $f(\ti\si)\in
\partial\h3$ is fixed by the elementary group $G_\si.$ We will
distinguish between the two cases when the group $\gsi$ is
loxodromic cyclic or parabolic subgroup of rank $2.$

\medskip

 \noindent  {\bf Case 1.  The group $\gsi$ is loxodromic.}

\medskip

 \noindent Let $V\subset\Pi$ be a regular neighborhood of the vertex $\si$. Then
the punctured neighborhood $V\setminus\si$ is homotopically
equivalent to the one-skeleton $L^{(1)}$ of the link $L$ of $\si$.

We will call {\it realization} of   $L$  a graph $\La\subset
V\setminus\si$ such that the canonical map $L\to \La$ is a
homeomorphism.   Let us fix
   a maximal tree $T$ in  $\La$, and let $y_i$ be the edges
from $\La\setminus T$ which generate the group $\pi_1(L)\
(i=1,...,k).$

By its very definition,  the $G$-equivariant  map $f:P\to\h3$
sends the edges of $P$ to geodesics of $\h3$. So let $\gsi=<g>$
and let $\ga\subset M$ be the corresponding closed geodesic in $
M$. We denote by $\ag\subset \h3$ the axis of the element $g$ and
by $g^{+}, g^{-}$ its fixed points on $\partial\h3.$ Let us assume
that $f(\ti\si)=g^{+}$. For   $X\subset M$ we denote by $\rm
diam(X)$ the diameter of $X$ in the hyperbolic metric of $M.$

Recall that the map   $f:P\to \h3\cup
\partial \h3$  constructed in Lemma \ref{emb} induces the map $F: \Pi' \to M$.
We start with the following:

\medskip

\noindent {\bf Step 1.} {\it For every $\eta >0$ there exists a
 realization $\La$ of $L$ in $\Pi$
 such that   for the maximal tree $T$ of $\La$ one has

 $${\rm diam}(F(T))<\eta,$$

\noindent Furthermore, for every edge $y_i\in \La\setminus T$ its
image $F(y_i)$ is contained in a $\eta$-neighborhood
$N_{\eta}(\ga)\subset M$ of the geodesic $\gamma$\ (i=1,...,k).}

\medskip

\noindent {\it Proof:} We fix a sufficiently small neighborhood
$V$ of a vertex $\si$ in $\Pi$ (the "smalleness" will be specified
later on). Let $\ti\si\in\nu^{-1}(\si)$ be its lift to $P$ and let
$\ti\La$ and $\ti T$ be the lifts of $\La$ and $T$ to a
neighborhood $\ti V\subset \nu^{-1}(V)$ of $\ti\si.$ We are going
first to show  that, up to decreasing $V,$  the image $f(\ti T)$
belongs to a sufficiently small horosphere in $\h3$ centered at
the point $g^+$.

Let   $\al$ be an edge of $\Pi$ having $\si$ as a vertex and
  $\ti\al$ be its lift starting at a point
$\ti\si.$   Then  $a=f(\ti\al)\subset\h3$ is the geodesic ray
ending at the point $\g+$, let $a(t)$ be its parametrization. For
a given $t_0$ we fix a  horosphere
 $S_{t_0} $   based at  $\g+$ and
passing through the point $a(t_0)$.  Suppose  there is a simplex
in $P$ having two edges $\ti\al=[\ti\si, s], \ti\al_1=[\ti\si,
s_1]$ at the vertex $\ti\si$ and an edge $[s, s_1]$ in $\La$. The
horosphere $S_{t_0} $ is the level set of the Busemann function
$\be_{\g+}$ based  at the point $\g+$. So for the geodesic rays
$a=f(\ti\al)$ and $a_1=f(\ti\al_1)$   issuing from the point $\g+$
we have that the points  $f(s)=a(t_0)$ and  $f(s_1)=a_1(t_0)$
belong to the horosphere $S_{t_0}.$ Proceeding in this way for all
simplices whose edges share the vertex $\si$,
   we obtain  that $f(\ti T^{(0)})\subset S_{t_0}\subset\h3.$ Since
   $\La$ is   finite, so is the tree $\ti T$. By choosing $t_0$ sufficiently
large    ($t_0 > \Delta$) we may assume that $d(\al_i(t_0),
\al_j(t_0)) <\eta$ and $d(\al_i(t_0), \ag) < \eta\ (i,j=1,...,k).$
We now connect all the vertices of $f(\ti T)$ by geodesic segments
$b_i\subset\h3$. By convexity, and up to increasing the parameter
$t_0$, we also have $ d(b_i, \ag)<\eta.$

By Lemma \ref{emb} the map $f$ sends  the lifts $\ti y_i\in \ti T$
of the edges $ y_i\in \La\setminus T$    simplicially to $b_i\
(i=1,...,k)$; and  $f$ maps $\gsi$-equivariantly  the preimage
$\ti\La=\nu^{-1}(\La)$ to $\h3.$ Hence the map $f$ projects to the
map $F : \La\to M$ satisfying the claim of Step 1. \bx

\bigskip

\noindent \noindent {\bf Step 2.} {\it Definition of the
polyhedron $\Pi\check\  $ }

\bigskip

Using the initial  orbihedron $\Pi$  we will construct a new
polyhedron $\Pi\check\ $ having the following properties :

\medskip

 a) $\Pi^{(0)}={\Pi\check\ }^{(0)}$  and
 $\Pi=\Pi\check\ $ outside of $V$;

\medskip

 b) $\pi_1(L^*)=\gsi,$ where $L^*$ is the link of $\si$ in $\Pi\check\ ;$

\medskip

 c) $\pi_1(\Pi\check\ )\cong G.$

\medskip

 \noindent  The graph  $\La$ realizes the link of the vertex $\si$ so
  there exists an epimorphism  $\pi_1(\La)\to  <g>.$  Every edge $y_i\in
\La\setminus T$ which is  a generator of the group $\pi_1\La$ is
mapped onto $g^{n_{y_i}}$ in $\gsi\ (i=1,...,k).$ We now subdivide
each edge $y_i$   by edges $y_{ij}\ (i=1,...,k, j=1,...,n_{y_i}),$
and denote by  $\La'$  the obtained graph. Let $S$ be a circle
considered as a graph with one edge $e$ and one vertex $u.$ Then
there exists a simplicial map from $\La'$ to $S$ mapping
simplicially each edge $y_{ij}$ onto $S$.

To construct polyhedron $\Pi\check\  $, we replace the
neighborhood $V$ by the cone of the above map. Namely, we   first
delete the vertex $\si$ from $\Pi$ as well as all edges connecting
$\si$ with $L$. Then we
  connect the vertices of the edge
  $y_{ij}$ with the vertex $u\in S$ by  edges which we  call
{\it vertical} $(i=1,...,k, j=1,...,n_{y_i})$. So $\Pi\check\ $ is
the union of $\Pi\setminus V$ and
   the rectangles $R_{ij}$, which are  bounded by   $y_{ij}$,
   two vertical
  edges   and the loop
  $S.$   The set of rectangles
  $\{R_{ij}\ \vert\ i=1,...,k,\ j=1,...,n_{y_i}\}$
  realizes the  epimorphism $\pi_1(L)\to G_\si.$   By Van-Kampen theorem
  we have $\pi_1(\Pi\check\  )\cong G,$
 and   the  conditions a)-c)
 follow. \bx

  \medskip

 \noindent {\bf Step 3.} {\it There exists a constant $c$
 (depending
 only on the topology of $\Pi$) such that for all $\eta >0$, there exists a
map $F\check\ :\Pi\check\ \to M$ such that

\begin{itemize}

\item[1)] $F\check\  $ induces an isomorphism on the fundamental
groups,

\medskip

\item [2)]   $F\check\  \vert_{\Pi\check\  \setminus V}=F$,

\medskip

\item [3)] $\displaystyle\sum_{ij} {\rm Area}(F\check\ (R_{ij})) <
c\cdot\eta.\hfill(2)$

\end{itemize}

}

 \bx

\medskip

\noindent {\it Proof:} We choose a neighborhood $V$ of the
singular point $\si$ and put $\displaystyle F\check\
=F\vert_{\Pi\setminus V}$. Using Step 2 we transform the
orbihedron $\Pi$ to $\Pi\check\  $ in the neighborhood $V$  and
let
 $P\check\ $ be the universal covering
 of $\Pi\check\  $. Note that, by construction, $P\check\ $ is obtained by adding
 the $G$-orbit of
 the rectangles $R_{ij}$ to the preimage  $\ti
 \La'=\nu^{-1}(\La')$ of the graph $\La'\ (i=1,...,k,\ j=1,...,n_{y_i}).$

 We will now
extend the map $f$  defined on $P\setminus V$  to the polyhedron
$P\check\ \setminus P$ as follows. We first subdivide every
segment $b_i$ in $n_{y_i}$ geodesic subsegments $b_{ij}\subset
b_i$ corresponding to the edges $y_{ij}$.  We now project
orthogonally each $b_{ij}$ to $\ag$  and let $\ti\ga\subset\ag$
denote its image. Let $\tau_{ij}\subset\h3$ be the rectangle
formed by $b_{ij}$, $\ti\ga$ and these two orthogonal segments
from $b_{ij}$ to $\ag$ whose lengths are by Step 1 less than
$\eta.$ We extend the map $f$ simplicially to a map $f\check\  $
sending  the rectangle $\nu^{-1}(R_{ij})$ to the rectangle
$\tau_{ij}\ (i=1,...,k,\ j=1,...,n_{y_i})$. Note that by
construction the lift $\ti S$ of the circle $S$ is mapped on $\ti
\ga.$ The map $f\check\ $ descends to a map $F\check\  :
\Pi_*\setminus \Pi\to N_{\eta}(\ga).$ It induces the epimorphism
$\pi_1\Pi\check\ \to G.$

Let us now make  the area estimates for the added rectangles
$\tau_{ij}$. Each rectangle $\tau=\tau_{ij}$ has four vertices $A,
B, C, D$ in $\h3$ where $B=gA, D=g(C)$ and the segment $[A,
B]\subset A_g$ is the orthogonal projection of $[C, D]$ on  $\ag$.
The rectangle $\tau$ is bounded by these two segments and two
perpendicular segments $l_1=[A, C]$ and $l_2=[B, D]$ to the
geodesic $\ag$ ($l_2=g(l_1)$). We have $\tau\subset ABC'D$ where
$\angle BDC'={\pi\over 2}$ and $\be=\angle BC'D < {\pi\over 2}$.
Then by [Be, Theorem 7.17.1] one has $\cos(\be )\leq \sinh (d(B,
D))\cdot \sinh l(\ga).$ Therefore ${\rm Area} (\tau ) < {\pi\over
2}-\be,$ and   $\sin ({\rm Area} (\tau ))\leq \sinh \eta\cdot
\sinh l(\ga).$ Summing up over all segments $b_{ij}$ we arrive to
the formula (2). This proves  Case 1. \bx

\bigskip

 \noindent  {\bf Case 2.  The group $\gsi$ is parabolic.}

\bigskip

 \noindent The proof is  similar   and even simpler   in this case.
 Let again $T$ be the maximal
tree of
   the graph $\La$ realizing the link
$L$ of the vertex $\si$. We start by   embedding  a lift $\ti
T^{(0)}$ of the zero-skeleton of $T^{0}$ into a horosphere
$S_{t_0}\subset\h3$ based  at the parabolic fixed point
$p\in\partial\h3$ of the group
 $\gsi=<g_1, g_2>\cong \Z+\Z.$  Then, using  Lemma \ref{emb},
  we construct an embedding   $f
:\ti \La^{(0)}\to S_{t_0}$ of the zero-skeleton of the graph $\ti
\La=\nu^{-1}(\La)$ into the same horosphere $S_{t_0}$ invariant
under $\gsi$ (which was not so in the previous case). Since the
number of vertices of $\ti T$ is finite, for any $\eta
>0$ we can choose a horosphere $S_{t_0}\ (t_0> \De)$  such that ${\rm
diam} \ti T < \eta$. Fixing a point $O\in S_{t_0},$  we can also
assume that $d(O,\ti T^{(0)})<\eta.$

 Now, let us modify the orbihedron $\Pi$ in the neighborhood $V$ of $\si.$
First we delete the vertex $\si$ from $\Pi$ and  all edges
connecting $\si$ with the graph $\La$.
 We then add to the obtained orbihedron   a torus
$\cal T$ with two intersecting loops $C_1$ and
  $C_2$ representing the generators of $\pi_1(T,u)$ where
    $u\in C_1\cap C_2$.   To realize the
  epimorphism $\pi_1\La\to \gsi$ in $M$ we proceed as before. For any edge
   $y\in \La\setminus T$
  corresponding to the element $g=ng_1+mg_2$ in $\gsi$ we add
  a rectangle $R$ bounded by $y$, two edges connecting  the end
  points of $y$ with $u$ and a loop $C\subset\T$ representing the element
  $g$ in $\pi_1(T,u)$. Let $\Pi\check\  $ denote the obtained
  orbihedron.

Coming back to $\h3,$ let us assume for simplicity that
  $p=\infty$ and the horosphere $S_{t_0}$ is   a Euclidean
  plane. By Lemma \ref{emb}
   the map $f$    sends the edges   $\ti y_i\in \ti\La\setminus \ti T$
   to the geodesic edges $b_i$ connecting the vertices
   of $f(\ti T).$

  We now construct the rectangles $\tau_i$ by projecting the
end points of the edges $b_i$ to the corresponding vertices of the
Euclidean lattice given by the orbit $\gsi O.$ Let us briefly
describe this procedure in case of one rectangle $\tau$. Suppose
that the edge $y\in \La\setminus T$ represents the element
$g=ng_1+mg_2\in\gsi.$ Let $A$ and $gA$ be vertices of $f(\ti T)$
belonging to $S_{t_0}$ connected by a geodesic segment $b$
corresponding to $y.$ Let $\tau\subset\h3$ be the geodesic
  bounded by the edges $b, l=[O, A], gl, gb.$
We extend the map $f\check\  :\ti R\to\tau$ where  $\ti R$ is a
lift of the corresponding rectangle  $R$  added to $\Pi.$ The map
$f\check\  $ descends now to a simplicial map $F\check\   :
\Pi\check\ \to M$ sending the torus $\cal T$ into a cusp
neighborhood of the manifold $M.$ Since the rectangle $\tau$
belongs to $\eta$-neighborhood of the horosphere $S_{t_0},$ its
area, being close to the Euclidean one, is bounded by
$c\cdot\eta^2$ for some constant $c>0.$ Summing up over all edges
$y_i$ we obtain that the area of added rectangles does not exceed
$k \cdot c\cdot \eta^2.$ This proves Case 2. \bx

\medskip

\noindent To finish the proof of Proposition \ref{polappox}, we
note that the initial orbihedron  $\Pi$ is finite, so it has a
finite number of vertices $v_1,...,v_l$ whose vertex groups are
either loxodromic or parabolic. So for a fixed $\ve >0,$ we apply
the above simplicial "blow-up" procedure in a neighborhood of each
vertex $v_i\ (i=1,...,l)$. Finally, we obtain a 2-complex
$\Si_\ve;$ and  the simplicial map $\phi_\ve  : \Si_\ve \to M$
which induces an isomorphism on the fundamental groups and such
that
 $\vert {\rm Area} (\varphi_{\ve}(\sie)) - {\rm Area}(f(\Pi'))\vert
< \psi(\eta),$ where $\psi$ is a continuous function such that
$\displaystyle\lim_{\eta\to 0} \psi(\eta)=0.$ So for $\eta$
sufficiently small we have  $\psi(\eta) <\ve$ which proves the
Proposition.\qed

\bigskip

\noindent {\it Proof of Theorem B.} Let $G$ be the fundamental
group of a hyperbolic 3-manifold $M$ of finite volume. Let
$\Pi=P/G$ be a finite orbihedron realizing the invariant $T(G,
E)$, i.e. $\pi_1^{\rm orb}(\Pi)\cong G$, all vertex groups of
$\Pi$ are elementary and $\vert\Pi^{(2)}\vert=T(G, E)$.  Hence
${\rm Area}(F(\Pi'))=\pi\cdot T(G,E).$ Then by Proposition
\ref{polappox} for any $\ve
>0 $ there exists a 2-polyhedron $\Si_{\ve}$ and a map $\psi_\ve
:\Si_\ve \to M$ which induces an isomorphism on the fundamental
groups and such that
$$ {\rm Area} (\psi_\ve(\Si_\ve)) <\pi T(G,E)+\varepsilon $$ By \cite{C}
we have ${\rm Vol} M < {\rm Area}(\psi_{\ve}(\sie))<\pi
T(G,E)+\varepsilon\ (\forall\ve >0).$ It follows ${\rm Vol} M \leq
\pi T(G,E).$ Theorem B is proved.\qed

\section{Proof of Theorem A.}

In this Section we finish the proof of
\bigskip

\bigskip

\noindent {\bf Theorem A.} {\it There exists a constant $C$
  such that  for every hyperbolic 3-manifold $M$ of finite volume
  the following inequality holds}:
 $$C^{-1}  T(G, E_{\mu})  \leq  {\rm Vol} (M) \leq  C T(G, E_{\mu})\hfill \eqno(*)$$

 \bx
\medskip

 The right-hand side of the inequality (*) follows from Theorem B
if one puts $E=E_\mu$. So we only need to prove the left-hand side
of (*). We start with the following Lemma dealing with
n-dimensional hyperbolic manifolds :

\begin{lem} \label{balls}  Let $M$ be a
n-dimensional hyperbolic manifold of finite volume. Then there
exists a 2-dimensional triangular complex $W\subset\Mtk$   such
that $\pi_1(W)\hookrightarrow\pimtk$ is an isomorphism and $$\vert
W^2\vert \leq \sigma\cdot {\rm Vol}(M),$$ where $\vert W^2\vert$
is the number of 2-simplices of $W$ and $\sigma=\sigma(\mu)$ is a
constant depending only on $\mu$.
\end{lem}

\proof The Lemma is a quite standard fact, proved for $n=3$ in
\cite{Th} and more generally in \cite{G}, \cite{BGLM}, \cite{Ge}.
We provide a short proof of it  for the sake of completeness.
Consider a maximal set of points $\A=\{a_i\ \vert\ a_i\in \Mtk, \
d(a_i, a_j]
> \mu/4\}$ where $d(\cdot,\cdot)$ is the hyperbolic distance
of $M$ restricted to $\Mtk.$  By the triangle inequality we obtain
$$B(a_i, \mu/8)\cap B(a_j, \mu/8)=\emptyset\ \ {\rm if}\ i\not=
j,$$ where $B(a_i, \mu)$ is an embedded ball in $M$ (isometric to
a ball in $\hn$)  centered at $a_i$ of radius $\mu.$ By the
maximality of $\A$ we  have  $\displaystyle\Mtk\subset {\cal
U}=\bigcup_i B(a_i, \mu/4).$ Recall that the nerve $N\U$ of the
covering $\U$ is constructed as follows. Let  $N\U^0=\A$ be the
  vertex set. The vertices $a_{i_1},...,a_{i_{k+1}}$ span a
$k$-simplex if for the corresponding balls we have
$\displaystyle\bigcap_{j=1}^{k+1} B(a_{i_j},
\mu/4)\not=\emptyset.$ Since the covering $\U$ is given by balls
embedded into $M$,   the nerve $N\U$ is homotopy equivalent to
$\U\ $ [Hat, Corollary 4G.3].

Note that $\Mtk\hook \U\hook\mtk.$ Indeed if $x\in \D B(a_i,
\mu/4)$ then by the triangle inequality we have  $B(x,
\mu/4)\subset B(a_i, \mu/2)$, and so both are embedded in $M.$
Then $x\in \mtk.$ By the Margulis lemma,  as the corresponding
components of their thin parts are homeomorphic, the embedding
$\Mtk\hook\mtk$ is a homotopy equivalence. It implies that the
complex $N\U$ is homotopy equivalent to $\Mtk.$ Let $W$ denote the
2-skeleton of $N\U.$ Then it is a standard topology fact that $W$
carries the  fundamental group of $N\U$ \cite{Hat}. Therefore,
$\pi_1 W\cong \pimtk$.

It remains to count the number of 2-faces of $W.$ We have for the
cardinality $\vert\A\vert$ of the set $\A$:

$$\vert\A\vert \leq  {\Vo(\Mtk)\over {\rm Vol}(B(\mu/8))}
\leq {\Vo(M)\over \Vo(B(\mu/8))}\ ,$$

\noindent where  $B(\mu)$ denotes a ball of radius $\mu$ in the
hyperbolic space $\hn$. The number of faces of $W$ containing a
point of $\A$ as a vertex is at most $\displaystyle
m={\Vo(B(\mu/2))\over \Vo(B(\mu/8))}.\ $ Then $$\displaystyle
\vert W^{(2)}\vert \leq C_m^2\ {\Vo (M)\over
\Vo(B(\mu/8))}=\sigma\cdot \Vo (M)\ ,$$
 where $\displaystyle \sigma=\sigma(\mu)=
{C_m^2\over\Vo(B(\mu/8))}$. This completes the proof of the Lemma.
\bx

\bigskip

  Suppose now that $M$ is a hyperbolic 3-manifold of
finite volume and let $\mu=\mu(3)$   be the 3-dimensional Margulis
constant. We are going to use a result of \cite{De} which we need
to adapt to our Definition \ref{relinv} of the invariant $T$. So
we start with the following:

\begin{rem}
\label{t0} In the definition of the invariant $T$ in \cite{De}
there is one more additional condition compared  to our Definition
\ref{relinv}. Namely, it requires that every element of a system
$E$ fixes a vertex of $P$. To be able   to use the results of
\cite{De} we will
  denote by $T_0(G, E)$ the invariant defined in \cite{De} and
keep the notation $T(G, E)$ for that of our Definition
\ref{relinv}. Notice that nothing  changes   for the absolute
invariant $T(G)$.
\end{rem}

Let
    $l_1,...,l_k$ be the
set of closed geodesics in $M$  of length  less than $\mu$. Then
by \cite{Ko} the manifold $\displaystyle M'=M\setminus \bigcup_i^k
l_i$ is  a complete hyperbolic manifold of finite volume and
$\pimtk\cong\pims$.

Let $\E_\mu $ denote the system $\pi_1(\partial \Mtk)$ of
fundamental groups of the boundary components of the thick part
$\Mtk$. We have the following :

\begin{lem}
\label{pres}

$$T_0(\pim, \E_\mu) \leq T_0(\pi_1(M'), \pi_1(\D M'))\leq T_0(\pim, \E_\mu) +2k.\eqno (5)$$

\end{lem}

 \proof 1) Consider first the left-hand side.
Let $G=\pim$ and $G'=\pi_1(M').$ Let $\E'_\mu=\{E_{k+1},...,E_n\}$
be the set of fundamental groups of cusps of $\Mth$. Let us fix a
two-dimensional $(G',\E'_\mu)$-orbihedron $P'$ containing $T_0(G',
\E'_\mu)$ triangular 2-faces. The pair $(G',\E'_\mu)$ acts on its
orbihedral universal cover $P'$ \cite{H}.  Let $N(l_i)$ be a
regular neighborhood of the geodesic $l_i\in M\ (i=1,...,k)$ and
$H_i=<\al_i, \be_i>$ be the fundamental group of the torus $T_i=\D
N(l_i)$ where $\al_i$ is freely homotopic to $l_i$ in $N(l_i).$
The group $H_i$ fixes a point $x_i\in P'.$ We will now construct a
2-orbihedron $P$ for the couple $(G,E_{\mu})$ as follows. The
group $G$ is the quotient of $G'$ by adding the relation $\be_i=1\
(i=1,...,k)$. We
 identify the vertices of $P'$ equivalent under  the groups
generated by $\be_i\ (i=1,...,k)$. The natural projection map
$P'\to P$ consists of contracting each  edge of $P'$ of the type
$(y, \beta_i(y))\ (y\in P'^{(0)})$ to a point. The projection  has
connected fibres so the 2-orbihedron $P$ is  simply connected and
the pair $(G,E_{\mu})$ acts on it. The procedure did not increase
the number of 2-faces, and we have : $\vert\Pi^{(2)}=P/G\vert\leq
\vert \Pi'^{(2)}=P'/G'\vert.$ Thus $T_0(\pim, E_{\mu})\leq
T_0(\pim', \pi_1(\D M')=\E'_\mu)$.

2)  Let $\Pi$ be the 2-orbihedron which realizes
$T_0(\pim,\E_{\mu}),$ and let $P$ be its universal cover. To
obtain a $(\pim', \E'_\mu)$-orbihedron
 we modify $P$ as follows. Let
$H_i=<h_i>$ be the loxodromic subgroup corresponding to the
geodesic $l_i\subset M$ of length less than $\mu\ (i=1,...,k).$
Let $x_i\in P$ be a vertex fixed by the subgroup $H_i$. Notice
that the group $G'$ is generated by $G$ and   elements $\be_i$
such that $[h_i, \be_i]=1\ (i=1,...,k)$. So we add to $\Pi$ a new
loop $\be_i$ (by identifying it with the corresponding element in
$G$) and  glue a disk whose boundary is the loop corresponding to
$[h_i, \be_i]$. By triangulating each such a disk  we add $2k$ new
triangles to $ \Pi^{(2)}$. Thus the universal cover $P'$ is
obtained by adding to $P$  a vertex $y_i$ and its orbit
$\{Gy_i\},$ so that the points $\be_ih_igy_i$ are identified with
$h_i\be_igy_i$. We further add the rectangle $gD_i\ (g\in G)$
whose vertices are $h_igy_i, \be_ih_igy_i, \be_igy_i, gy_i$ and
  subdivide it by  one
of the diagonal edges, say $(h_igy_i,\be_igy_i)\ (i=1,...,k)$. The
construction gives a new 2-complex $P'$ on which the pair $(G',
\E'_\mu)$ acts simplicially. We claim that $P'$ is simply
connected. Indeed if $\al$ is a loop on it, since $P$ is simply
connected, $\al$ is homotopic to  a product of loops belonging to
the disks $gD_i$ so $\al$ is a trivial loop. Since the
2-orbihedron $\Pi'=P'/G'$ contains $\vert\Pi^{(2)}\vert+2k$ faces,
we obtain $T_0(\pim', \pi_1 (\D M'))\leq T_0(\pim, \E_{\mu})+2k$
which was promised. \qed

\begin{rem}
It is worth   pointing out that in the context of volumes of
hyperbolic 3-manifolds the following  inequality (similar to (5))
is known:

$$ {\rm Vol} (M) < {\rm Vol} (M') < k\cdot (C_1(R)\cdot {\rm Vol}(M)
+C_2(R)),\hfill\eqno(\dagger)$$

\noindent where $R$ is the maximum of radii of the embedded tubes
around the short geodesics $l_i\ (i=1,...,k)$ and $C_i(R)$ are
functions of $R\ (i=1,2).$ The left-hand side of ($\dagger$) is
classical and due to W.~Thurston \cite{Th},  the right-hand side
is proved recently by I.~Agol, P. A.~Storm, and W.~Thurston
\cite{AST}\bx
\end{rem}

\bigskip

\noindent {\it Proof of the left-hand side of the inequality (*):}
By Lemma \ref{balls} the thick part $\Mtk$ of $M$ contains a
2-dimensional complex $W$ such that $\pi_1 W\hook\pimtk$ is an
isomorphism and $\vert W^{(2)}\vert < \sigma\cdot {\rm Vol} (M) $
for some uniform constant $\sigma.$ Consider now the double
$N=D\Mtk$ of the manifold  $\Mtk$ along the boundary $\D\Mtk.$  By
repeating the argument of Lemma \ref{balls} to each half of $N$ we
obtain two complexes $W$ and $\tau(W)$ embedded in $N$ where $\tau
: N\to N$ is the involution such that $\Mtk=N/\tau.$ By Van-Kampen
theorem the fundamental group of the complex $V=W\cup \tau(W)$ is
generated by $\pi_1W$ and $\pi_1(\tau(W))$ and  is isomorphic to
$\pi_1(N)$. Furthermore, for the number of two-dimensional faces
in $V$ we have $\vert N^{(2)}\vert= 2\vert W^{(2)}\vert$. So by
Lemma \ref{balls} $T(\pi_1 N)\leq \vert V^{(2)}\vert <
2\sigma\cdot {\rm Vol} (M).$ The group $\pi_1 N$ splits as the
graph of groups whose two vertex
 groups are $\pimtk$. The edge groups of the graph of groups
are given by the system $\E_{\mu}$. As $\pimtk\cong\pims$ and $M'$
is a complete hyperbolic
   3-manifold of finite volume it follows from  Lemma \ref{red}
   that the above splitting is reduced and rigid.
So  by \cite{De} we have:

$$T(\pi_1 N) \geq 2 T_0(\pimtk,\E_\mu).\eqno(6)$$

\noindent Then by Lemma \ref{pres} $T_0(\pimtk,\E_\mu) \geq
T_0(\pim, \E_\mu),$ and therefore $$\si^{-1}\cdot T_0(\pi_1(M),
\E_\mu) < {\rm Vol} (M).$$

\noindent Recall that the initial system $E_\mu$ of elementary
subgroups includes all elementary subgroups of $\pi_1(M)$ whose
translation length is less than $\mu.$ So $\E_\mu\subset E_\mu$
implying that $T(\pi_1(M), E_\mu) \leq T_0(\pi_1(M), \E_\mu)$. We
finally obtain

$$C^{-1}\cdot T (\pi_1(M), E_\mu) < {\rm Vol} (M),$$

\noindent where $\displaystyle C= \sigma.$ The left-hand side of
(*) is now proved. Theorem A follows. \bx

\section{Concluding remarks and questions.}

The finiteness theorem of Wang affirms that there are only
finitely many hyperbolic manifolds of dimension greater than $3$
having the volume bounded by a  fixed constant \cite{W}. So it is
natural to compare the volume of a hyperbolic manifold
$M=\hn/\Gamma$ with the absolute invariant $T(\G).$ In the case
$n>3$ the
 inequality $${\rm const}\cdot T(\G) \leq {\rm Vol} (M)$$   follows from [Ge,
 Thm 1.7] (see also  Section 2 above, where instead of  $T(\pim, E)$
 one needs to consider $T(\pi_1 (M))$ and use the
 fact that $\pimtk\cong \pi_1(M)$).
 However, the result \cite{C}
 is not known in higher dimensions. Thus we have the following :

 \begin{ques}
 \label{ineq}

Is there a constant $C_n$ such that for every lattice $\G$ in
$\iso$   one has

 $${\rm Vol} (\Gamma) \leq C_n\cdot T(\G)\ ?$$

\end{ques}

\begin{rem} (M.~Gromov) The answer is positive if $M$ is a compact
hyperbolic manifold of dimension 4. Indeed in this case by the
Gauss-Bonnet formula one has ${\rm Vol} (M) = {\Omega_4\over
2}\cdot \chi(M),$ where $\Omega_4$ is the volume of the standard
unit 4-sphere. Hence ${\rm Vol} (M) < {\Omega_4\over 2}\cdot (2-
2b_1 + b_2)$ where $b_i={\rm rank}(H_i(M, \Z))$ is the $i$-th
Betti number of $M\ (i=1,2).$ Since $b_2 < T(\pi_1(M)),$ one has
${\rm Vol} (M) < {\Omega_4\over 2}\cdot (2+b_2) < \Omega_4\cdot
T(\pi_1(M))\ ({\rm as}\ T(\pi_1(M))
> 1).$

\end{rem}

\noindent Recently it was shown by D.~Gabai, R.~Meyerhoff, and
P.~Milley that the Matveev-Weeks 3-manifold $M_0$ is the unique
closed 3-manifold of the smallest volume \cite{GMM}. Furthermore,
C.~Cao and R.~Meyerhoff found cusped 3-manifolds $m003$ and $m004$
of the smallest volume \cite{CM}, \cite{GMM}. In this context we
have the following :

 \begin{ques}

Is the    invariant $T(\pi_1(M), E_{\mu})$ on the set of
 compact hyperbolic 3-manifolds attained on the manifold $M_0$  ?
 Is the minimal relative invariant
$T(\pim, E_{\mu})$ on the set of cusped finite volume 3-manifolds
  attained on the manifolds $m003$ and $m004$ ?

\end{ques}

 \noindent Thomas Delzant: IRMA, Université de Strasbourg, 7 rue R. Descartes , 67084
Strasbourg Cedex, France; delzant@math.u-strasbg.fr

\bigskip

\noindent Leonid Potyagailo: UFR de Math\'ematiques, Universit\'e de Lille 1,
59655 Villeneuve d'Ascq cedex, France; potyag@gat.univ-lille1.fr

\clearpage

\end{document}